\documentclass{amsart}

\usepackage{amsfonts,amssymb,amsmath,amsxtra,graphicx}
\usepackage[mathscr]{eucal}

\newtheorem{thm}{Theorem}[section]

\theoremstyle{remark}

\theoremstyle{plain}

\newcommand*{\SideSet}[2]{%
  \ensuremath{%
    \hphantom{^{#1}}%
    {\vphantom{#2}}^{\mathpalette\myllap{#1}}%
    #2%
  }%
}%
\newcommand*{\myllap}[2]{\llap{$#1#2$}}%

\def\brep{\mathbf b} % band representation
\def\Bd{\operatorname{Bd}} % boundary of a manifold
\def\Int{\operatorname{Int}} % interior of a manifold
\def\sw{\operatorname{sw}}  % self-winding (Laufer)
\def\GCD{\operatorname{GCD}} 
\def\gr{\operatorname{graph}}
\def\id#1{\operatorname{id}_{#1}}
\def\Im{\operatorname{Im}}
\def\pr{\operatorname{pr}}   
\def\Re{\operatorname{Re}}
\def\Reg{\mathscr R}   
\def\Sing{\mathscr S}   
\def\sub{\subset} 
\def\transp#1#2{(#1\thickspace#2)} % transposition in permutation group
\def\unclose#1{#1\spcheck}
\def\C{\mathbf C} \def\P{\mathbf P} \def\R{\mathbf R} \def\Z{\mathbf Z} 
\def\Ss{\mathcal S} 
\def\b{\beta} 
\def\bhat{\widehat\b}
\def\d{\delta}
\def\e{\varepsilon} 
\def\g{\gamma} 
\def\th{\theta}
\def\s{\sigma} 
 
\def\D{\Delta}
\def\G{\Gamma}
\def\le{\leqslant}
\def\ge{\geqslant}

\thanks{Research partially supported by NSF Grant MCS 76-08230.
This survey was originally published in 
\emph{L'Enseignement Math\'ematique 29} (1983), 185--208, 
and in \emph{N{oe}uds, Tresses, et Singularit\'es}
(Monographie No. 31 de L'Enseignement Math\'ematique), 
ed.\ C. Weber (Kundig, Geneva, 1983), 99--122.
The present {\LaTeX}ed redaction, completed in June 2001,
corrects several typographical errors without, I hope, introducing new 
ones; some minor emendations and references to post-1983 results 
have been added as footnotes, but I have not made a thoroughgoing 
effort to update the text.  Anyone
having further information---for instance, on the current status of 
questions incorrectly described as still open---is encouraged to let 
me know by e-mail to {\tt lrudolph\@black.clarku.edu}.}

\begin{document}
\title{Some knot theory of complex plane curves}
\author{Lee Rudolph}

\maketitle

\section{Aspects of the ``placement problem'' \\
for complex plane curves}

How can a complex curve be placed in a complex surface?

The question is vague; many different ways to make it more specific
may be imagined.  The theory of deformations of complex structure,
and their associated moduli spaces, is one way.  Differential geometry
and function theory, curvatures and currents, could be brought in.
Even the generalized Nevanlinna theory of value distribution, for 
analytic curves, can somehow be construed as an aspect of 
the ``placement problem''.

By ``knot theory'' I mean to connote those aspects of the situation
that are more immediately topological.  I hope to show that there is
something of interest there.

\section{A triptych.}

Here are three ways to interpret the phrase ``knot theory of complex
plane curves''.

Globally: the ``complex plane'' is projective space $\C\P^2$
or affine space $\C^2$; a ``curve'' is an algebraic curve
(in projective space) or an algebraic or \-analytic curve (in
affine space); here, ``knot theory'' has historically been
largely 
concerned with studying the ``knot group'', though there 
are also results on ``knot type''.

Locally: a ``complex plane curve'' is the germ of a plane curve
(algebraic, analytic, or formal) over $\C$; this is the study of
singularities, and ``knot theory'' has been the classical knot
theory of links in the $3$-sphere, put to work in the service
of that study.

In between: a ``complex plane curve'' is an analytic curve in 
a reasonable open set in a complex surface (chiefly, in the 
theory as so far developed, the interior of a ball or a bidisk),
well-behaved at the boundary; a knot-theorist can study either
of two codimension-$2$ situations---the complex curve in its
ambient space, or the boundary of this pair.

This middle panel of the triptych has been less studied than
the other two, though it is of obvious relevance to both.

\section{R\'esum\'e of basic definitions}

By \emph{complex surface} I mean a smooth manifold of $4$ real dimensions, 
equipped with a complex structure.  A \emph{complex curve} $\G$ in a 
complex surface $M$ is a closed subset which is locally of the form 
$\{(z, w) \in U \sub \C^2:f(z, w) = 0\}$ where $f: U \to \C$ is 
a nonconstant complex analytic function.  A \emph{Riemann surface} is a 
smooth manifold of $2$ real dimensions, equipped with a complex 
structure.

It is a fundamental fact, to which is due the especial 
appositeness of classical knot theory to the study of curves in 
surfaces, that any complex curve $\G\sub M$ has a \emph{resolution} 
of the following sort: There is a Riemann surface $R$, and a 
holomorphic mapping $r: R \to M$, so that $r(R) = \G$; in fact, 
there is a discrete (possibly empty) subset $\Sing(\G) \sub \G$,
the \emph{singular locus} of $\G$ in $M$, so that the \emph{regular locus}
$\Reg(\G) = \G\setminus \Sing(\G)$ is a Riemann surface, and $R$ is the union 
(with what turns out to be a unique topology and complex 
structure) of $\Reg(\G)$, on which $r$ is the identity, and a discrete 
set $r^{-1}(\Sing(\G)) \sub R$ mapping finitely-to-one onto 
$\Sing(\G)$.

The singular locus is, of course, exactly the set of points of $\G$ 
at which, no matter what the local representation of $\G$ as the 
zeroes of an analytic function $f(z,w)$, the (complex) gradient 
vector $\nabla f$ vanishes.

If $P$ is a point of $\G$, and $Q \in r^{-1}(P) \sub R$, then the 
germ at $P$ of the \hbox{$r$-image} of a small 
disk on $R$ centered at $Q$ is called a \emph{branch} 
of $\G$ at $P$.  (Abusively, ``branch'' may also be 
used below to refer to some representatives of this germ.) 
Naturally, at a regular point there is only one branch; but there 
may be only one branch at a point, and the point still be 
singular.

References: \cite{G-R}, \cite{Mi1}.

\section{Local knot theory in brief }

Using local coordinates in the resolution $R$ and the ambient 
surface $M$, one sees that each branch of a curve $\G$ can be 
parametrized either by $z = t, w = 0$ or (more interestingly) by some 
pair $z = t^m, w = t^n+c_{n+1} t^{n+1}+\dots c_N t^N$,
$t\in\C$, with $n > m$.  (In the original choice of coordinates, $r$ 
might well involve genuine power series; but it is not hard to make a 
formal change of coordinates to one of the forms above, involving 
only polynomials; and it is not much harder to prove a comparison 
theorem, the remote ancestor of that of M.\  Artin, which shows 
that actually the formal change of coordinates can be taken to be 
somewhere convergent.) Consider the ``approximations'' to such a 
branch, gotten by dropping all terms of $w$ from some degree on up: 
so the first approximation is $(t^m,t^n)$, and the $(N-n+1)$st 
is the branch we began with.  Each of these is itself a map onto a 
branch of some curve, generally not one-to-one.

Define integers $g(1), \dots, g(N-n+1)$ by saying that the 
$k$-th approximation is $g(k)$-to-one in a punctured neighborhood of 
$t=0$.  Then $g(1)=\GCD(m,n)$, $g(k+1)$ divides $g(k)$, and 
$g(N-n+1) = 1$.  These integers can be calculated as follows.  
Let $\C[[t]]$ be the algebra of formal power series, with unique 
maximal ideal $m = t\C[[t]]$.  
Let $A_k$ be the $m$-adically closed subalgebra 
generated by $1$ and the components of the $k$-th approximation.  
Then $g(k)$ is the least integer $g$ such that 
$A_k\sub\C[[t^g]]\sub \C[[t]]$.  (One gets the same answer 
starting from the algebra $\C\{t\}$ of somewhere-convergent 
power series.) A parametrization of the branch covered by the 
$k$th approximation is 
$z=t^{m/g(k)}, w=t^{n/g(k)}+ \dots +c_{n+k-1}t^{(n+k)/g(k)}$.

The knots in which we are interested arise when we intersect the 
branch 
under investigation with the boundary of an infinitesimal 
$4$-disk containing the singular point.  The $4$-disk used may be 
either a \emph{round disk} $D_\e^4 = \{(z, w): |z|^2 
+|w|^2 = \e^2\}$ with boundary the \emph{round sphere} $S^3_\e$, 
or a \emph{bidisk} $D(\e_1, \e_2) 
= \{(z, w): |z| \le \e_1, | w | \le \e_2\}$, 
with boundary comprised of two solid tori 
$\Bd_1 D(\e_1, \e_2) = \{|z|=\e_1,|w|\le\e_2\}$
and $\Bd_2 D(\e_1, \e_2) = \{|z|\le\e_1, |w| = \e_2\}$ 
which together make up a $3$-sphere with corners.  Whether one 
uses round disks or bidisks, one obtains a knot of the same type.  
The bidisk boundary is more convenient here, when we are studying 
the branch parametrically; from the 
assumption that $n > m$ we can 
see that, for sufficiently small $\e > 0$, the branch 
intersects $\Bd D(\e,\e)$ only along $\Bd_1 D(\e, \e)$.

The first approximation to the branch actually meets 
$\Bd_1 D$ on the torus $\{|z|=\e ,|w|=\e^{n/m}\}$, 
where it covers, $g(1)$ to one, a \emph{torus knot} 
of type $O\{m/g(1),n/g(1)\}$.  
(Here is the notation I am using, cf.\ \cite{Ru4}.  
If $K$ is any \hbox{oriented} knot in an oriented $3$-sphere, 
with closed tubular neighborhood $N(K)$, let $L$ be an 
oriented simple closed curve on $\Bd N(K)$ which is not 
null-homologous on this torus; then there are relatively 
prime integers $p$ and $q$ so that $L$ has linking number 
$q$ with $K$ and represents $p$ times the class of $K$ in $H_1(N(K);\Z)$.  
We then call $L$ a \emph{cable of type $(p,q)$ about $K$} 
and denote it by $K\{p, q\}$.  When cabling is
iterated, excess curly braces become semicolons.  The unknot is 
denoted by $O$; a cable about the unknot is also called a 
\emph{torus knot}; a cable about \dots a cable about 
the unknot 
is an \emph{iterated torus knot}.) This knot type does not 
change when $\e$ is made smaller.

Now suppose that for all sufficiently small $\e > 0$, the $k$-th 
approximation to a branch intersects $\Bd D(\e, \e)$ in a knot of 
type $O\{p_1,q_1;\dots ; p_k,q_k\}$.  
Considering how we pass to 
the next approximation we see that there are relatively prime 
integers $p_{k+1}$ and $q_{k+1}$ so that, for all sufficiently small 
$\e>0$, the $(k+1)$-st approximation to the branch intersects 
$\Bd D(\e, \e)$ in a knot of type 
$O\{p_1, q_1 ;\dots ; p_k, q_k; p_{k+1}, q_{k+1}\}$.  
(The difference between successive approximations is 
$0$ or a monomial $c_{n+k}t^{n+k} \not\equiv 0$, which contributes 
an ``epicycle'' that for small enough $\e$ precisely creates a cabling.)
In fact, $p_{k+1}= g(k)/g(k+1)$ (note that for any $K$ and $q$, 
$K\{1, q\}$ is the same knot type as $K$); the formula for 
$q_{k+1}$ is more complicated, and we won't give it 

Consider a curve with a singular point at which there are 
two or more branches.  Coordinates in the ambient surface can be 
chosen so that each branch differs only by a diagonal linear 
transformation in $(z,w)$ from one of the form just studied 
(including the non-singular case $z=t, w=0$).  Each 
branch individually contributes an iterated torus knot to the 
\emph{link of the singularity} $\G\cap\Bd D(\e, \e)$; and in fact
they all fit together nicely.  An elegant description of how 
they do is given in \cite{E-N-preprint}\footnote{Published 
as \cite{E-N-published}.}; 
see also, and for this section genera1ly, 
\cite{Le} and \cite{Mi2} and references cited therein.

After Burau, Zariski, \emph{et al.}, had established that 
any point of a curve in a (nonsingular) surface had local 
topology that was completely described by this link-type 
invariant, the strictly topological investigation of singular 
points seems to have languished for some decades.  (The algebraic 
geometers, of course, had also established that this link-type 
invariant---more precisely, the sequences of pairs $(p_i, q_i)$
for each branch, and the linking numbers between the 
iterated torus knots of different branches from which the 
whole link of the singularity can be reconstructed---was 
equivarient to some numerical invariants which had long 
been known and which could be detected purely algebraically,
namely, the Puiseux pairs of the various branches and the 
intersection multiplicity of the pairs of branches.  They also
pressed forward with their investigations of continuous invariants 
within the family of singularities of a given link type.  
But that is another story.) However, in the late 1960's,
Milnor \cite{Mi2} gave new life to the subject when 
he showed that the link of a singularity was a 
``fibred'', or Neuwirth--Stallings, link.

Milnor's proof uses the round-sphere model.  He shows that, if 
$\G\sub\C^2$ is the zero-locus of $p(z,w)\in\C[z,w]$, $p(0,0)=0$, 
then for all sufficiently small $\e>0$, the restriction $\phi$ of 
the map $\arg p:\C^2\setminus \G \to S^1:(z, w) \mapsto p(z,w)/|p(z, w)|$ 
to $S^3_\e\setminus \G$ is the projection map of a fibration over $S^1$.
The fibre is diffeomorphic to the interior of the surface 
$F_0 = S^3_\e \cap \{(z, w):p(z, w) \text{ is real and non-negative}\}$.
(Note that the change in viewpoint from bidisk boundary to round sphere 
is accompanied by a change from branch-as-parametrized-disk to 
branch-as-level-set.)

We will see below that the link of a singularity is in a natural 
way a closed strictly positive braid; I will give a geometric 
proof of the well-known fact that such a closed braid is a fibred 
link.

Inspired by Milnor's Fibration Theorem, a number of 
mathematicians began investigations of knot-theoretical 
properties of the links of singularities.  The
fibration $\phi$ determines an autodiffeomorphism 
of $F_0$ (fixed on the boundary),
unique up to isotopy relative to the boundary, which 
is variously called the \emph{characteristic map}, \emph{holonomy}
or \emph{monodromy} of the fibration; it induces an
automorphism (also called the monodromy) of the integral homology 
of $F_0$.  From the homology monodromy one can calculate the Alexander 
polynomial of the link of the singularity; this was done in \cite{Le},
where it was also shown that two branches defined iterated torus 
knots in the same knot-cobordism class if and only if they defined 
knots of the same knot type, the proof following from a
study of the roots of the Alexander polynomials.

I wondered how independent these distinct knot-cobordism classes 
might be, in the knot-cobordism group; in particular, I asked 
 \cite{Ru6} whether the equation $[K_0]=\sum\limits_{i=1}^n [K_i]$,
in which $[K_i]$ represents the (non-trivial) knot-cobordism
class of the link of a singular branch, $i = 0,\dots,n$ had any 
solutions other than $K_1=K_0$, $n=1$.  Litherland, using 
his calculations of the signatures of iterated
torus knots \cite{Li}, was able to show that there were only such 
trivial solutions.  It follows that, for instance, there is 
no family $\{\G_s\}$, $|s|< \e$ , of (local) curves in a
small hall in $\C^2$ so that $\G_s$ for $s\ne 0 $ has two 
singular points each with a single branch while $\G_0$ has 
only one singularity, locally of the form $z = t^2, w = t^5$.  
Is there another proof of the non-existence of such a deformation? 
(Multiplicities would allow two cusps.)

Litherland's formulas, of course, give all the various 
signatures of the links of singularities (though the expression 
is in closed form only by the use of a counting function involving 
``greatest integer in\dots'', which makes them rather a
bore to calculate).  If one lowers one's sights, and asks only about 
the classical signature (that corresponding to the root $-1$ of unity), 
and then only about its sign, an easy direct proof---again, using the 
representation of the link as a closed 
positive braid---shows that the signature of the link of a 
singularity is positive, \cite{Ru5}.

Finally, some conjectures on less algebraic knot invariants of 
links of singularities should be mentioned.  The \emph{Milnor 
number} $\mu$ of a singularity is the
rank of $H_1(F_0 \Z)$.  Let us look at a single branch, for 
convenience.  Then Milnor
conjectured \cite{Mi2} that $\mu/2$, which is the genus of 
$F_0$ and therefore (by a general
theorem about fibred links) the genus of the knot $\Bd F_0$,
actually is the \emph{slice genus} of $\Bd F_0$.  One can make the 
weaker conjecture that at least $\mu/2$ is the \emph{ribbon genus} of
$\Bd F_0$.  Milnor also wondered if this integer equalled the 
\emph{\"Uberschneidungszahl}, or \emph{gordian number}, of $\Bd F_0$;
again the conjecture can be weakened if one introduces 
the concepts of ``slice \"Uberschneidungszahl'' and ``ribbon
\"Uberschneidungszahl'', cf.\ \cite{Ru2}.  The conjectures are 
true in various cases where direct calculations can be made
(e.g., the cusps $z=t^2, w = t^3)$, but I know of no general 
results.\footnote{Kronheimer and Mrowka, 
by proving the local Thom Conjecture \cite{K-M}, answered Milnor's
question affirmatively.  See also \cite{Ru7} for further
knot-theoretical consequences of the truth of the local
Thom Conjecture.}

\section{Global knot theory in brief---the projective case}

A curve $\G \sub \C\P^2$ can be given by its resolution 
$r: R \to \G$ (a complex-analytic map from a compact 
Riemann surface into $\C\P^2$ which is generically 
one-to-one on $R$) or by its polynomial 
$F(z_0,z_1,z_2)\in \C[z_0,z_1,z_2]$ (the 
homogeneous polynomial of least degree, not identically 
zero, which vanishes at every point of $\G$).  
These suggest different kinds of knot-theoretical 
questions.  One can consider all curves with diffeomorphic 
resolutions (the requirement that the curves have complex-
analytically equivalent resolutions would be too stringent, and 
is less topological), and ask how differently they can be placed 
in the plane.  Or one can consider families of curves, each cut 
out by a polynomial of some fixed degree.
Let $P_d$ denote the projective space of the vector space of 
homogeneous complex polynomials in $(z_0, z_1, z_2)$ of degree $d$.
Because we never want to consider curves with multiple 
components, we throw out of $P_d$ the algebraic subset 
corresponding to reducible polynomials with a multiple factor;
the remaining Zariski-open subset $Q_d$ corresponds to the set of 
what we may call curves of \emph{geometric degree} $d$.  
If (the equivalence class of) $F(z_0, z_1, z_2)$ belongs to
$P_d$, let 
$\G_F = \{(z_0:z_1:z_2) \in \C\P^2:F(z_0, z_1, z_2) = 0)$; 
then $F \in Q_d$ if and only if there is an open dense set of 
lines in $\C\P^2$ which intersect $\G_F$ transversely in $d$ distinct 
points.

The condition that $\G_F$ have a singular point is, of course, an 
algebraic condition on $F$.  Let $S_d \sub P_d$ be the algebraic 
subset of singular curves without multiple components, and $R_d = 
Q_d\setminus S_d$ the Zariski-open subset 
of polynomials of \emph{geometrically 
regular curves} of geometric degree $d$.  Any curve $\G_F \in R_d$ 
is its own resolution ($r$ = identity).  By connecting any two 
$F, G \in R_d$ with a path in $R_d$, one may construct an isotopy (which 
may be effected by an ambient isotopy) between the curves $\G_F$ 
and $\G_G$ in $\C\P^2$; so all these curves are diffeomorphic, and 
of the same knot type in the plane.  More generally, $F \in Q_d$ lies 
in a maximal connected subset of $Q_d$ of polynomials $G$ such that 
$\G_F$ and $\G_G$ are ambient isotopic, through algebraic curves.  
These subsets form a stratification of $Q_d$ which is little 
understood.  Zariski \cite{Z} showed that two (singular) curves in 
$Q_6$, homeomorphic and with the same type and number of 
singularities (cusps), were not in the same stratum, by showing 
that the knot groups $\pi_1(\C\P^2\setminus \G_F)$ and 
$\pi_1(\C\P^2\setminus \G_G)$ were not isomorphic.  
In general, as we will see below, the knot group cannot 
distinguish strata.

An interesting question (I do not know to whom it is due: I heard 
of it in Dennis Sullivan's problem seminar at M.I.T.\  in the 
summer of 1974) is whether there are curves $\G_F$ and $\G_G$ which 
are ambient isotopic but not so through algebraic curves.  I know 
of no results here.

The incidence structure of this stratification of $Q_d$ by 
``algebraic ambient isotopy types'' is, especially, not understood: 
this is the theory of degenerations.  It can be proved that the 
knot group associated to a given stratum is the homomorphic image 
of the knot group associated to any stratum incident to the given 
stratum.  Partly, it was the desire to apply this fact to the 
proof of the Zariski Conjecture (see below) which led 
investigators for many years to the study of some particular 
(unions of) strata to which we now turn.

First we recall the two simplest sorts of singularities.  A cusp 
has a sing1e branch, locally given by $z = t^2, w = t^3$; the link of 
a cusp is a trefoil knot (of a fixed handedness once one 
establishes conventions).  A node has two branches, each itself 
nonsingular, with distinct tangent lines; it can be locally given 
by the equation $zw = 0$, and its link is a Hopf link of two 
components (linking number $+1$).  A curve $\G$ is a \emph{node curve}
if all its singularities (if any) are nodes, and a \emph{cusp curve}
if all its singularities are either nodes or cusps.

We also recall, what we have not needed before, the notion of 
reducibility: a curve $\G$ is \emph{reducible} if its resolution is not 
connected; alternatively $\G_F$ is reducible if and only if $F$ is 
reducible but square-free.  A curve that is not reducible is 
\emph{irreducible}.

The extreme of reducibility is displayed by any $\G\in Q_d$ which is 
the product of $d$ linear factors.  Then the curve $\G$ is the union 
of $d$ projective lines, which we will say (here) are in \emph{general 
position} if $\G_F$ is a node curve, that is, if no three of the 
lines are concurrent.  Let $L_d\sub Q_d$ be the set of all such 
comp1ete1y reducible curves.  Then $L_d$ is a single stratum.  Let 
$N_d \sub Q_d$ be the set of polynomials of node curves; $N_d$ is a 
union of strata.  What is now called the Severi 
Conjecture\footnote{As of 1986, a theorem of Harris \cite{Ha}.}
is the 
statement that $L_d$ is incident to every stratum in $N_d$; in other 
words, that every node curve can be degenerated to $d$ lines in 
general position.  We will compute the knot group of $d$ lines in 
general position below.  It is, in particular, abelian.  
Consequently, the truth of the Severi Conjecture would imply that 
\emph{the knot group of any node curve is abelian}---a statement long 
known as the Zariski Conjecture, which has recently been proved 
true by quite other means \cite{F-H, De}\footnote{An alternative
proof of a stronger theorem (the link at infinity may be any
closed positive braid) was given in 1988 by Orevkov
\cite{O}, using braid-theoretic methods related to those in 
\cite{Ru1}.}.  Of course, independent 
of the truth of the Severi Conjecture, one can study the union 
$M_d \sub N_d$ of those strata which actually are incident to $L_d$.  
Moishezon \cite{Mo} calls $M_d$ the \emph{mainstream} of node curves 
in his investigation of ``normal forms for braid monodromies''.
Such normal forms (when they exist) enrich the datum of the knot 
group by giving it in a particularly nice presentation related to 
the algebraic geometry.

Now let $K_d \sub Q_d$ correspond to the cusp curves.  Here the knot 
groups need no longer be abelian.  In fact, for 
$$
F(z_0, z_1, z_2) = z_1^2 z_2^2+4z_0(z_2^3-z_1^3) 
+ 6z_0^2 z_1 z_2-27z_0^4
$$
in $K_4$, a curve with three cusps and no nodes (which has 
resolution 
$$
r: \C\P^1 \to \G_F: (t_0: t_1; t_2)\mapsto 
(t_0^2 t_1^2: t_1^4+2t_0^3 t_1:2 t_0 t_1^3+t_0^4)),
$$
the knot group can be computed (as by Zariski \cite{Z} or, 
algebraically, by Abhyankar \cite{Ab}) to have the presentation 
$(a, b:aba = bab, a^4 = 1, a^2 = b^2)$, making it non-abelian of 
order $12$.

The knot groups of cusp curves have been studied because of their 
application to the study and possible classification of complex 
(algebraic) surfaces.  In fact, if $f:Y\to \C\P^2$ is a so-called 
stable finite morphism, $\Sigma' \sub Y$ the locus where $f$ is not 
\'etale, $\Sigma = f(\Sigma')$, then $\Sigma$ is a cusp curve.

Zariski commissioned van~Kampen, in the early 1930's, to 
calculate the knot group of an arbitrary curve \cite{vK}; van
Kampen gave his solution in terms of a certain presentation of the 
knot group.  If $\G$ has (geometric) degree $d$, then 
van Kampen's presentation has $d$ generators $x_1,\dots,x_d$
which represent loops in a fixed projective line $\C\P^1_\infty$
transverse to $\G$; the intersection $\G \cap \C\P^1_\infty$ contains 
$d$ points $P_1, \dots, P_d$, and $x_i$ is a loop from a basepoint 
$*\in \C\P^1_\infty$ out to $P_i$, around it once counterclockwise, 
and back to $*$.  One relation is then that $x_1\dotsm x_d = 1$.  
The rest arise by carrying $\C\P^1_\infty$ around certain loops of 
lines.  In fact, let $\C\P^{2*}$ be the dual projective plane, 
each point of which is a line in $\C\P^2$; and let $\G^*$ contain 
all lines which are either tangent to $\G$ or pass through one of 
its singular points.  Then $\G^*$ is a curve in $\C\P^{2*}$.
If $*$ and $\C\P^1_\infty$ are sufficiently general, then the 
pencil of lines in $\C\P^{2}$ through $*$, which is itself a line in 
$\C\P^{2*}$, will be transverse to $\G^*$.  The (free) fundamental 
group of the complement of $\G^*$ in this pencil is naturally 
represented in the automorphism group of the free group 
$(x_1, \dots, x_d:x_1\dotsm x_d = 1)$.  The 
rest of the relations needed for the van~Kampen presentation 
of $\pi_1(\C\P^{2}-\G; *)$ come, then, by declaring this 
representation trivial.  One obtains a finite presentation,
of course, by choosing generators of the acting free group.
Moishezon's problem of ``normal forms'' is essentially the 
problem of making a good choice.  Several modernizations \cite{Abe}, 
\cite{Che}, \cite{Cha} of van~Kampen's proof have been 
published in recent years.

In a \emph{standard} van~Kampen presentation (where the generators of 
the acting free group are free generators), each relation 
corresponds either to a singularity of $\G$ or to a simple vertical 
tangent to $\G$; and (up to the action of the corresponding free 
generator) each relation is of a certain canonical form, which 
depends only on the \emph{closed braid type} (\S7) of the link of the 
branch(es) at the point of $\G$ through which the line in the 
pencil passes that gives the relation in question, where this 
line itself is used to find the axis of the closed braid.  In 
particular, the knot group of a node curve always has a standard 
van~Kampen presentation in which each relation either sets 
conjugates of two $x_i$ equal (from a simple vertical tangent) or 
says that two such conjugates commute (from a node); if
``conjugates'' could be deleted, the Zariski Conjecture would be 
trivially true.

There is also a great body of work on ``knot groups'' of curves in 
(compact, smooth complex surfaces other than $\C\P^{2*}$ and on 
the related issue of fundamental groups of surfaces; we cannot 
touch on these topics here.

\section{Global knot theory in brief---the affine case}

Little appears to be known about algebraic curves in affine 
space, from the knot-theoretical viewpoint.  The gross algebraic 
topology (even just homology theory) of $\C\P^{2}$ is implicated 
with the quite rigid geometry; but affine space is contractible, 
and on the other hand its geometry is ``infinite'' (for instance in 
the sense that there are Lie groups of arbitrarily high dimension 
contained in the group of biregular automorphisms of $\C^2$), so 
that the conspirators have fallen out and neither can give away 
much about the other.

One might think, for example, to study the embedding of a curve 
$\G$ in $\C^2$ by first embedding $\C^2$ itself into $\C\P^{2}$. 
Then the affine complement $\C^2\setminus\G$ becomes the projective 
complement $\C\P^{2}\setminus (\G \cup \C\P^1_\infty)$, 
where $\G \cup \C\P^1_\infty$ is a (reducible) projective algebraic curve.  
The obstacle to this program is the unfortunate fact that $\C^2$,
just as an algebraic surface, without distinguished coordinates, 
is not uniquely embedded as $\C\P^2\setminus \C\P^1_\infty$ .  Any 
biregular automorphism of $\C^2$ (in particular, one of the vast 
majority which cannot be extended regularly to $\C\P^{2}$) will 
move $\G$ around, and so the configuration of $\G \cup \C\P^1_\infty$ 
is not determined by the embedding of $\G$ in $\C^2$.  (For 
instance, though the geometric number of points at infinity 
on $\G$ is determined by $\G$, the algebraic intersection number 
of the closure of $\G$ with the line at infinity can be made 
arbitrarily large.  Likewise the local singularities at infinity 
are not determined by the affine curve.)

The main theorems known here have been proved by Abhyankar and 
his collaborators \cite{A-M, A-S}\footnote{See also \cite{Suz}, 
for an analytic proof of the main theorem of \cite{A-M}, published 
slightly earlier.}.  They are unknotting theorems, 
in the sense that they take this form: ``Let $\G$ be a certain 
curve in $\C^2$, and let $i: \G \to \C^2$ be any algebraic embedding; 
then there is a biregular automorphism of $\C^2$ returning $i$ to 
the inclusion map''.  Briefly, such a curve $\G$ cannot be knotted 
in $\C^2$.  

However, for most of the curves they deal with, these theorems 
are not genuinely topological, for the reimbedding $i$ is required to 
be an embedding \emph{of $\G$ with its given structure as a variety},
and generally there might be moduli.  Only in the original theorem 
\cite{A-M} (which had been stated, but not correctly proved, by
Segre) are there no conceivable moduli, when $\G$ is a straight 
line.  Then the theorem is this.

\begin{thm}
Let $\G \sub \C^2$ be an algebraic curve without singularities, 
homeomorphic to $\C$.  Then there is a biregular change of 
coordinates $A:\C^2 \to \C^2$ so that $A\G$ is a straight 
(complex) line.  
\end{thm}

A topological proof has been given in \cite{Ru4}.\footnote{At
least two other topological proofs have since been given.
That in \cite{N-R} dispenses with knot cobordism and uses 
instead a notion of ``unfolding a fibred knot''; 
that in \cite{N} uses the calculus of splice diagrams.
Both these proofs have the further virtue that they recover, 
not just the Abhyankar--Moh---Suzuki classification of 
polynomial embeddings of $\C$ in $\C^2$, but also the 
Za{\u\i}denberg--Lin classification of singular polynomial
injections of $\C$ in $\C^2$, \cite{Z-L}.} It goes like this.  
One shows (just as for a singular point) that the intersection of 
$\G$ (which we can assume to be parametrized by $z = p(t)$, $w = q(t)$, 
$p, q \in\C[t]$) with a \emph{very large} bidisk boundary is an iterated 
torus knot $K = O\{m_1, n_1 ; \dots ; m_s, n_s\}$, with 
$m_1 = m/\GCD(m,n), n_1 = n/\GCD(m, n), m = \deg p, n = \deg q$.  
By hypothesis, $K$ is a slice knot.  This forces $K = O$, in 
particular, one of $m_1, n_1$ is 1.  Thereafter the argument 
is as in \cite{A-M}---if (say) $m_1 = 1$ and $p$ and $q$ are 
monic then the biregular change of coordinates 
$(z, w) \mapsto (z, w-z^{m/n})$ carries $\G$ to 
another curve satisfying the hypotheses, of lower bidegree; and 
so we proceed until one of $z, w$ is linear and the other constant.

As to analytic curves in affine space, almost nothing is known.
The obvious analogue of the Theorem above is definitely false: 
for it is known that the unit disk in $\C$ can be properly 
analytically embedded in $\C^2$ \cite{H}; since the disk and the 
line are analytically inequivalent, no analytic change of 
coordinates in $\C^2$ could unknot the disk to a line.  It is, 
however, perfectly possible that every such disk is smoothly 
unknotted.  Presently I am unable even to prove that an analytic 
line in $\C^2$ is smoothly unknotted.

\section{The middle range }

We return, as at the beginning of \S4, to the study of 
intersections of curves in $\C^2$ with round disks $D^4_r$ and 
their boundaries $S^3_r$, and bidisks $D(r_1, r_2)$ and their
boundaries.  Now the (bi)radii are no longer required to be very 
small.

An embedding $i:(S, \Bd S)\to (D^4_r, S^3_r)$ of a 
surface-with-boundary $S$ into a round disk is a 
\emph{ribbon embedding} provided that $N\circ i$ is a Morse 
function without local maxima on $\Int S$, where 
$N(z, w)=|z|^2+|w|^2$; and a surface-with-boundary 
$S \in D^4_r$, with $\Bd S = S^3_r \cap S$, 
is a \emph{ribbon surface} if the inclusion 
$(S, \Bd S)\sub (D^4_r, S^3_r)$ is isotopic through 
embeddings of pairs to a ribbon embedding.  
To demand that a surface be ribbon is to place genuine topological 
restrictions on the embedding.

A theorem of Milnor \cite{Mi1}, specialized to our dimensions, 
shows that if $\G\sub \C^2$ is a nonsingular analytic curve 
then for almost all choices of origin and radius, the inclusion 
of $(\G\cap D^4_r, \G\cap S^3_r)$ into $(D^4_r, S^3_r)$ is 
a ribbon embedding.  A continuity argument easily shows that 
for no matter what choice of origin, $N|\G$ has critical points, 
possibly degenerate, of index no greater than $1$.  It is easy 
to see that if $\G$ has singularities, an analogous theorem holds 
for $N \circ r: R \to [0, \infty{[}$ on the resolution.  All 
these results generalize the Maximum Modulus Principle.  Nothing 
much more seems to be known about big round disks and complex 
plane curves\footnote{Our ignorance is now much less extensive.
See footnote~2 and, especially, \cite{B-O}.}.

Turning our attention to bidisks, we let the way that they 
separate the variables $z$ and $w$ suggest an attitude to adopt 
towards our curves: consider one variable (conventionally $w$) as 
an analytic but possibly multiple-valued function of the other.

More precisely, let $E_n$ be the space of unordered $n$-tuples of 
points of $\C$ (duplications allowed).  Then $E_n$ inherits a 
topology, and a structure of algebraic variety (affine, and 
singular if $n\ge 2$), from its description as $\C^n/\Ss_n$, 
where the symmetric group $\Ss_n$ acts by permuting coordinates.  
Let $E_n$ keep its topology, but normalize and resolve its 
algebraic variety structure, by using the map $\C^n\to E^n$ 
which carries $(c_1,\dots, c_n)$ to $\{r_1, ..., r_n\}$ such that 
$(w-r_1)\dots (w-r_n) = w^n+c_1 w^{n-1} +\dots+c_n$.  
Now any function $F: X\to E_n$ can be called \emph{an $n$-valued 
(complex) function on $X$}.  The \emph{graph} of an 
$n$-valued function on $X$ is 
the obvious subset of $X \times \C$; adjectives like 
``continuous'', ``analytic'', ``algebraic'' apply to $n$-valued 
functions in the obvious way.

We make the convention that (if $X$ is not discrete) the entire 
image $F(X)$ should not lie in the subset $\D \sub E_n$ of 
unordered $n$-tuples with at least one duplication.
$\D$ is an algebraic hypersurface (irreducible, and 
singular if $n \ge 3$) in the affine space $E^n$, 
called the \emph{discriminant locus}. Its complement $E^n\setminus \D$ is 
called the \emph{configuration space (of $n$ distinct 
points in $\C$)}.

To allow infinity as a value, we could 
replace $\C$ by $\C\P^1$, $E_n$ by $\C\P^n$, and so on.

Let $f(z, w) \equiv f_0(z)w^n+f_1(z)w^{n-1}+\dots+f_n(z)\in\C[z,w]$.  
Historically \cite{Bl} the equation $f(z, w) = 0$ (or equivalently 
the curve it defines) was said to give $w$ as an algebraic function 
of $z$, provided only that $f(z, w)$ was without repeated factors and
without factors of the form $z-c$.  (Also, of course, $f_0(z) 
\not\equiv 0$.) Then, in fact, on the 
complement in $\C$ of the zero-locus of $f_0(z)$, the assignment 
$z \mapsto \{w :f(z, w) = 0)$ is an algebraic $n$-valued 
complex function.  A zero of $f_0(z)$ is called a pole of the 
algebraic function, and can be accounted for by letting infinity 
be a value.

If $f_0(z)\equiv 1$, so that there are no poles at all, the algebraic 
function is \emph{entire}.  More generally, if $f_0(z), \dots, f_n(z)$ 
are allowed to be entire functions of $z$ (in the usual
sense), then $f(z,w)=0$ gives $w$ as an $n$-valued \emph{meromorphic 
function}; and if also $f_0(z) \equiv 1$, $w$ is an \emph{entire 
analytic $n$-valued function}.  The graph of an $n$-valued
entire function is a curve (algebraic or analytic as the case 
may be); when there are poles the graph must be closed up to 
provide fibres over them.

Conversely, any algebraic curve in $\C^2$ becomes, after almost 
any linear change of coordinates, such a graph for some $n$.  
(This is not so for analytic curves, 
in general.) Thus we can 
study plane curves by studying certain curves in $E_n$.

Let $\g \sub \C$ be a simple closed curve, $R$ 
the compact simply-connected region it bounds, 
$F:R \to E_n$ a continuous $n$-valued function 
analytic on $\Int R$ with $F(\g) \cap \D = \emptyset$.  
Then there is some radius $M>0$ so that the graph of 
$F|\g$ lies in the open solid torus $\g\times\{w\in\C:|w|< M\}$; 
and this graph is a (not necessarily connected) $n$-sheeted 
covering space of $\g$. An application of one version of the 
Maximum Modulus Principle \cite{G-R} shows that actually the 
graph of $F$ itself is contained in 
$D = R\times\{w:|w|\le M\}\sub\C^2$, 
a topological $4$-ball (with boundary $3$-sphere piecewise as smooth as
$\g$).  Now, $F(R)\cap\D$ must be finite; 
let $F^{-1}(\D)\sub R$ be called the \emph{branch locus},
and denoted $B$.  One can easily see that the graph of $F$ 
in $D$ is a $2$-dimensional pseudomanifold-with-boundary 
(i.e., geometric relative cycle), with any singularities 
lying in $B\times\{w:|w|\le M\}\sub\Int D$; its boundary 
in $\Bd D$ is exactly the link $L$ which is the graph of
$F|\g$.  Furthermore, the graph of $F$ is naturally oriented 
(by its complex structure at the regular points), so $L$ has 
a natural orientation, and the projection $L\to\g$ preserves
orientations.

At this point it is convenient to introduce braids; a general 
reference is \cite{Bi}. The \emph{braid group on $n$ strings} 
is the fundamental group $B_n =\pi_1(E_n-\D; *)$ of the configuration 
space.  Let $l:[0, 2\pi] \to E_n\setminus \D$, $l(0) = l(2\pi)$, be a 
parametrization of a loop in the configuration space.  Then the 
graph of $l$ in $[0, 2\pi]\times \C$ is a \emph{geometric braid}, 
that is, the union of disjoint arcs, on which $\pr_1$ is a covering 
projection to $[0, 2\pi]$, and such that the unordered $n$-tuples 
of top and bottom endpoints are identical; each arc is called a 
string.  Under the map $[0, 2\pi ]\times\C\to S^1\times\C:
(\th, w)\to (e^{i\th}, w)$, a geometric braid is carried to a 
\emph{closed braid} in the open solid torus.  When $S^1\times\C$ is 
identified with the tubular neighborhood of an unknotted circle 
in $S^3$, in such a way that distinct circles 
$S^1\times \{z_0\}$ and $S^1\times\{z_1\}$ are (algebraically, and 
therefore geometrically) unlinked, then any closed braid becomes 
a knot or link in $S^3$, and it is naturally oriented.  For 
$\b\in B_n$, any closed braid constructed in this way from a 
loop which represents $\b$ is denoted $\bhat$.  
If, conversely, $L\sub S^1\times \C$ is an oriented link 
on which $\pr_1$ is an orientation-preserving $n$-sheeted 
covering map, then any choice of a basepoint $e^{i\th}\in S^1$ 
yields a loop in $E_n\setminus\D$, based at 
$* =\{w \in\C:(e^{i\th}, w)\in L\}$, 
and thus a braid $\unclose L\in B_n =\pi_1(E_n-\D;*)$ 
with $(\unclose L)\sphat = L$.

Since $\D$ is irreducible, the abelianization of $B_n$ is infinite 
cyclic, and in fact $B_n$ is normally generated by one element, 
that is, generated by a single conjugacy class.  Choose for the 
basepoint $*$ of $E_n -\D$ the (real) $n$-tuple $\{(1,\dots, n\}$.  
Let 
$$
g_i(z, w)\equiv\left(w^2-(2i+1)w+\left(i^2+i+\frac14(1-z\right)\right)%
\cdot\prod_{\underset{j\ne i, i+1}\to{j=1}}^{n-1}(w-j)\in\C[z,w],
$$
for $i=1,\dots,n-1$; and let $G_ i :\C \to E_n$  be the $n$-valued 
function corresponding to $g_i(z,w)$. 
If $R=\{z: |z|\le 1\}$, then each $G_i|R$ is an 
embedding of $R$ as a normal disk to $\D$ (at a regular point), 
with center 
$$
G_i(0)=\left\{1,\dots,i-1,i+\frac12,i+\frac12, i+2,\dots,n\right\}
$$
on $\D$, and basepoint $G_i(1)=*$. Giving $\Bd R$ its positive 
(counterclockwise) orientation, we get oriented loops in $E_n-\D$, 
and the homotopy class of $G_i(\Bd R)$ is denoted by $\s_i$ and 
called \emph{the $i$-th standard generator of $B_n$}.  (The geometric 
braids corresponding to the given construction are the standard 
pictures of the $\s_i$.) The set of standard generators does, in 
fact, generate $B_n$, cf.\ \cite{Bi}.  Each $\s_i$ is conjugate to 
$\s_1$.  Following \cite{Ru2}, let any braid in $B_n$ 
conjugate to $\s_1$ be called a \emph{positive band in $B_n$}; 
a loop in the configuration space represents a positive 
band if and only if it is the oriented boundary of an oriented 
disk in $E_n $ which meets the discriminant locus transversely in a 
single positive (regular) point.  The inverse of a positive band 
is a \emph{negative band}.

An ordered $k$-tuple $\brep = (b(1),\dots, b(k))$ of bands in 
$B_n$ is a \emph{band representation of length $k$ of the 
braid $\b(\brep) = b(1)\dotsm b(l)$}.  (A braid word is 
a band representation where each band is a standard generator 
or the inverse of a standard generator.)  Each braid has many band 
representations, corresponding to the various null-homotopies,
transverse to $\D$, of a loop representing the braid in $E_n -\D$
to a point in $E_n$.  (See \cite{Ru2} for a precise statement and 
proof.) Such a null-homotopy gives a map of a disk into $E_n$, 
transverse to $\D$---the length of any corresponding band 
representation is the geometric number of intersections of the 
disk with $\D$, and the number of positive (resp., negative) bands 
is the number of positive (resp., negative) intersections with 
$\D$.  In particular, suppose each such intersection is positive, 
so each band $b(s)$ is positive.  Then $\brep$, $\b(\brep)$, 
and the closed braid $\bhat(\brep)$ are all called 
(in \cite{Ru1,Ru2,Ru3,Ru4}) \emph{quasipositive}.  
The closed braid $L$, associated to an analytic $n$-valued 
function $F$ and a simple closed curve $\g$ 
which bounds a simply-connected in the domain 
of $F$, is quasipositive.  (If $F$ as given is not transverse to 
$\D$ in $R$, almost any small translation of $F$ in $E_n$ will 
become so while the braid type of $L$ won't change; 
and complex analytic intersections are positive.)

Conversely, it is shown in \cite{Ru1} that for every 
quasipositive band representation $\brep$ in $B_n$, there are 
an algebraic $n$-valued function and simple closed curve yielding 
the given band representation in the manner just exposed.  It is 
also shown (and this is why we have excluded poles) that any type 
of closed braid whatever can occur as the graph over $S^1$ of a 
meromorphic (algebraic) $n$-valued function on $\C$.  (But note 
that when poles actually do occur inside the simple closed curve, 
the closed braid is never the complete boundary of the piece of 
analytic curve inside a bidisk; a typical example is given by 
$f (z, w)\equiv zw-\frac14$, in $D(1, 1)$, $\g = S^1$.)

Let $e:B_n\to\Z$ be abelianization.  Thus $e(\b)$ is the 
\emph{exponent sum} of $\b$, when $\b$ is written as 
a braid word in the standard generators; or more generally it 
is the number of positive bands in $\brep$ minus the number 
of negative bands in $\brep$ when $\b(\brep) = \b$.  
Geometrically, $e(\b)$ is the linking number of (any loop 
representing) $\b$ with $\D$, in $E_n$.  
Analytically, $e(\b)$ can be given by an integral formula, as by 
Laufer \cite{Lau}, where it is called \emph{self-winding} (and is 
generalized to links that aren't necessarily given as closed 
braids).

When $\brep$ is quasipositive, $e(\beta(\brep))$ 
is the length of $\brep$, a fact with the following 
geometric meaning.  When $F:R \to E_n$ is smooth and 
transverse to $\D$, then the graph of $F$ is a smooth 
surface in $R\times\C$; the intersections with 
$\D$ correspond to ``simple vertical tangents'' to the graph, 
and projection from the graph of $F$ back to $R$ is a branched 
covering, with only two sheets coming together over each branch 
point in $R$.  Thus the Euler characteristic 
$\chi(\gr F)$ equals $n\chi(R)-l$, if $l$ 
is the number of branch points.  When $R$ is a disk and 
$F$ corresponds to a quasipositive band representation $\brep$ 
then $l$ is the length of $\brep$ and we recover a genus formula 
for the graph of $F$ in terms of $n$, the number of components of 
the boundary of the graph, and the exponent sum of the 
boundary.  More generally, when $F$ is analytic, even if it is not 
transverse to $\D$ it will have a well-defined positive 
intersection multiplicity at each point of intersection, which 
will equal the number of geometric intersections of almost any 
small (analytic) perturbation of $F$; thus its graph, which will 
now be a singular curve, will have well-defined multiplicities 
for each singular point, and again a genus formula can be 
recovered, this time involving also these multiplicities: 
cf.\ \cite{Lau}.

A very interesting subclass of the quasipositive braids consists 
of the \emph{positive braids}.  A braid in $B_n$ is positive if it can be 
written as a word in the standard generators without using their 
inverses, \emph{strictly} positive if each of $\s_1, \dots \s_{n-1}$ 
actually occurs.  Positive braids play an important algebraic role 
in the braid group (cf.\ \cite{Bi}).  Closed positive braids enjoy 
various nice knot-theoretical properties (cf.\ \cite{St}, \cite{Ru5}), 
and have turned up in diverse contexts---as knotted orbits 
of some special dynamical systems \cite{Bi-W}; and, what is 
relevant here, as the links of singular points of plane curves.

\begin{figure}
\centering
\includegraphics[width=4.5in]{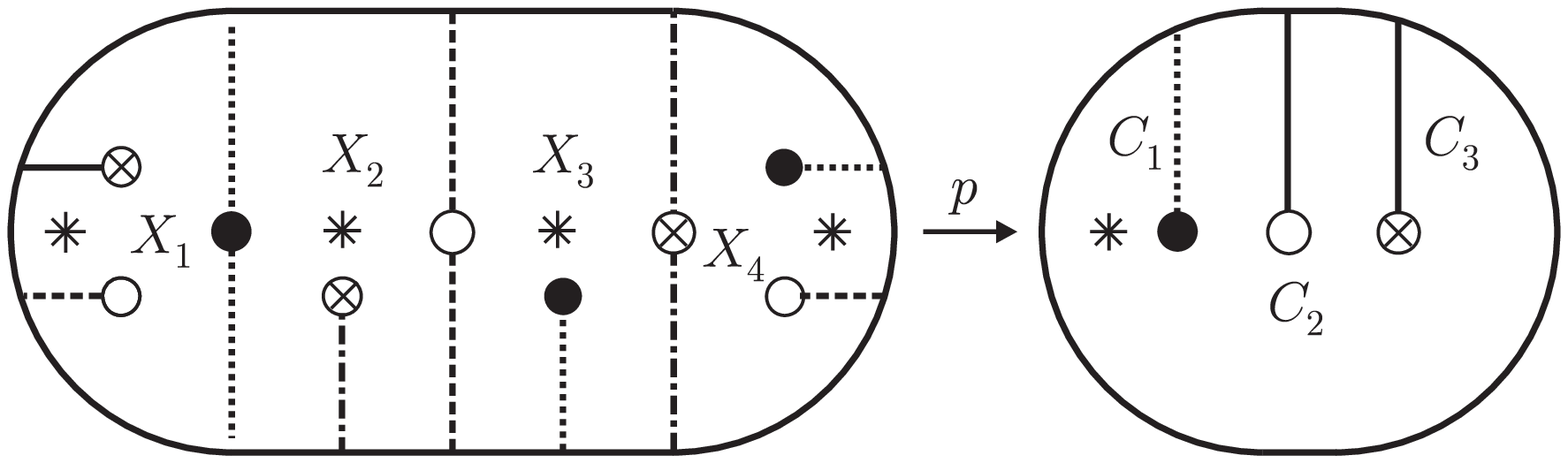}
\caption{$(n = 4)$}
\label{Figure 1}
\end{figure}

Let $f(z, w)\in\C[z,w]$ be squarefree, not divisible by $z$, 
and satisfy $f(0,0)=0$.  Then for $\e>0$ suffciently small, 
$f(z,w)=0$ defines an $n$-valued analytic function 
$F :\{z:|z|\le\e\}\to E_n$ with $F^{-1}(\D)=\{0\}$.  
Let $w_1(z),\dots,w_n(z)$ be the $n$ numbers in $F(z)$; 
then it is readily seen that the assignment 
$z\mapsto\{w_i(z)-w_j(z):1\le i,j\le n, i\ne j\}$ is an 
$n(n-1)$-valued analytic function.  
Without loss of generality, 
we may take $n$ and $\e$ so that $w_1(0)=\dots=w_n(0)=0$, 
and $w_i(z)-w_j(z)\ne 0$ for $z\ne 0$, $|z|\le\e$.  Now a 
straightforward calculation shows that for $z\ne 0$, $|z|\le\e$,
we have $d(\arg(w_i-w_j))/d(\arg z)>0$.  Consider the closed 
braid $L$, which is the graph of 
$F|\{z:|z|=\e\}$, and the 
link of the singularity of $\{f=0\}$ at $(0,0)$.  A braid diagram 
for $L$ may be obtained by projecting its ambient solid torus 
$S^1\times\C$ onto $S^1\times e^{i\th}\R$ orthogonally; for 
almost all $\th$ this will be a braid diagram in general position, 
from which a braid word may be read off in the usual way; and the 
signs of the crossings are precisely determined as the signs at 
the appropriate points of $d(\arg(w_i-w_j))/d\th$.  Since
$\th =\arg z$, the link of a singularity is a positive closed braid.  
In fact, it can be seen to be strictly positive, for if it were 
not, it would be a split link, in particular it would have 
components with zero algebraic linking---but the linking number 
of two components of the link of a singularity is the 
intersection number of the corresponding branches and is 
strictly positive.
It is known that a strictly positive closed braid is a fibred 
link, cf.\ \cite{St}, \cite{Bi-W}, which provides another proof
(in this dimension) of Milnor's Fibration Theorem (that the link 
of a singularity is fibred---Milnor, of course, gives an actual 
analytic formula for the fibration).  Here is a simple proof which
geometrically constructs a fibration of the complement of a 
strictly positive closed braid.  Let $p:X\to\C$ be the 
$n$-sheeted branched covering with branch locus 
$\{1,\dots,n-1\}$, where the permutation at $j$ 
is the transposition $\transp j{j+1}$.  Then $X$ is 
homeomorphic to $\C$ again.  For concreteness we realize 
$p$ as in Figure~\ref{Figure 1}: cuts
$C_j=\{w:\Re w=j, \Im w \ge0\}$ are made in the base space;
we coordinate $X$ so that the singular point of $p^{- 1}(j)$ is 
$j$, and so that $\{z :\Re z=j\}$ is one component of 
$p^{-1}(C_j)$; then the components of 
$p^{-1}(\C\setminus\cup_{j=1}^{n-1} C_j)$ 
are the sets 
$$
X_1=\{z: \Re z < 1\}, 
X_2 = \{z: 1 <\Re z < 2\},\dots, 
X_n = \{z: n-1 <\Re z\},
$$
known in the classical style as \emph{sheets} of the branched cover.  

\begin{figure}
\centering
\includegraphics[width=4.5in]{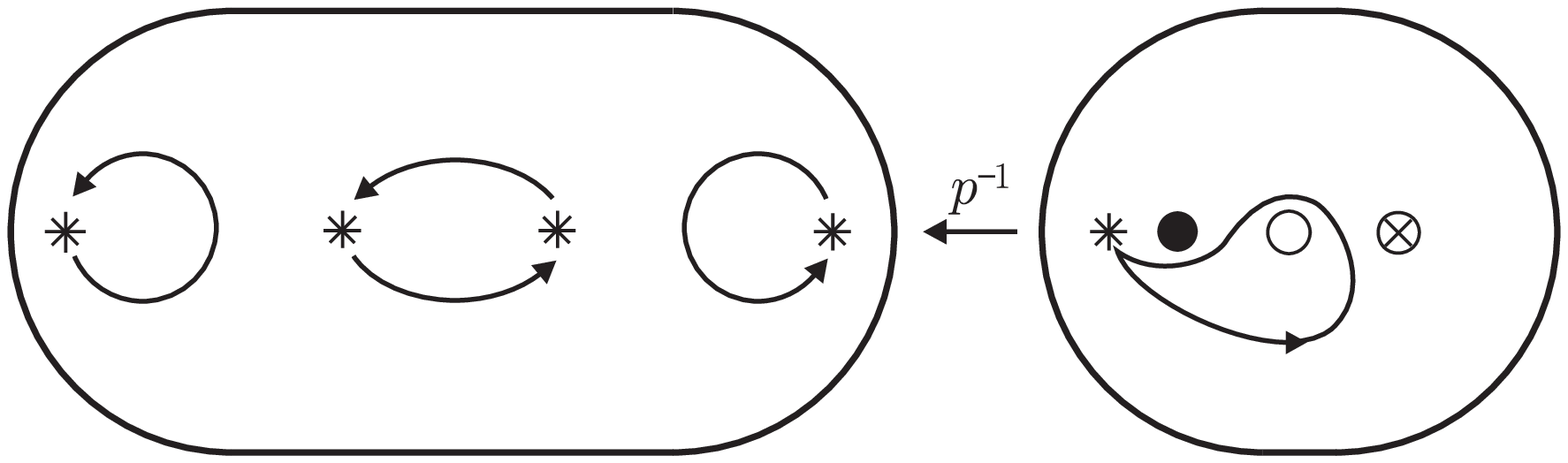}
\caption{$(n = 4)$}
\centerline{$(p^{-1})_*(x_2)=\s_2$}
\label{Figure 2}
\end{figure}

\begin{figure}
\centering
\includegraphics[width=4.5in]{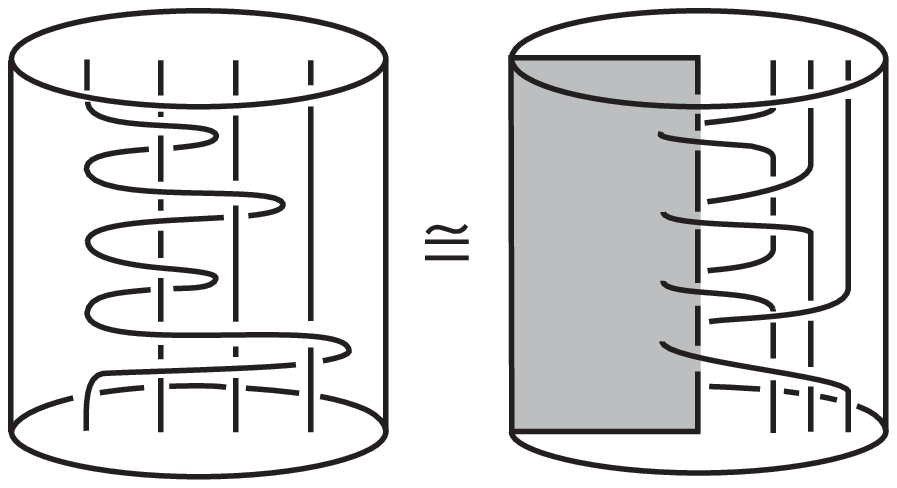}
\caption{$(n = 4)$}
\centerline{$v=x_1 x_2 x_1 x_3$}
\label{Figure 3}
\end{figure}

Now if we consider $E_n-\D$ to be the configuration space of $X$, the 
inverse of the covering map defines a continuous function from 
$\C-\{1,\dots,n-1\}$ into $E_n-\D$, inducing a homomorphism from 
the free group $\pi_1(\C\setminus\{1,\dots, n-1\}; 0)$ to the braid 
group $\pi_1(E_n-\D;p^{-1}(0))$.  One readily checks that this 
homomorphism is onto, carrying the obvious free generator $x_j$ 
of the free group (Figure~\ref{Figure 2}) to 
the standard generator $\s_j\in B_n$.  
Let $v=x_{j(1)}\dotsm x_{j(k)}$ be any strictly 
positive word in $x_1,\dots,x_{n-1}$,
$$
\b=\s_{j-1}\dotsm \s_{j(k)}=(p^{-1})_{*}(v)
$$
its strictly positive image in $B_n$.  We use $v$ to construct an 
auxiliary closed braid in $S^1\times\C$, the closure of 
$v'=A_{1,j(1)}\dotsm A_{1,j(k)}\in B_{n+1}$, 
where 
$$
A_{1,j}=(\s_1\dotsm\s_{j-1})\s_j^2(\s_1\dotsm\s_{j -1})^{-1}
$$
is one of the standard generators $A_{i,j}$ of the 
\emph{pure braid group} (cf.\ \cite{Bi} or see below).  
Now, $v'$ be realized as a geometric braid in two 
special ways: the first string can be made to wind in and out 
among the others, which are all straight; or the first string may 
be made straight, while the others wind around it in a succession 
of loops (Figure~\ref{Figure 3}).  On the first interpretation,
identifying the straight strings with 
$[0, 2\pi]\times\{1, ..., n-1\}$, the 
winding first string becomes the graph of a loop
$$
l: ([0,2\pi],\{0,2\pi\})\to(\C-\{1,\dots,n-1\}, 0)
$$
in the homotopy class $v$; and its inverse image under the 
branched covering 
$\id{S^1}\times p:S^1\times X \to S^1\times\C$ is 
a geometric braid representing $\b$.  On the second
interpretation, identifying the single straight string 
with $[0,2\pi]\times\{0\}$, and
taking care that each other string winds monotonically around 
this axis, the fibration of $S^1\times (\C-\{0\})$
over $S^1$ by $(e^{i\th},w)\mapsto\arg w$ lifts back through the
branched covering to a fibration of $(S^1\times X)-\bhat$ over $S^1$.  
(The strictness is used at this point, to ensure that in fact 
there is a non-zero winding number for each string.
Positivity, however, could be weakened to ``homogeneity''
in the sense of \cite{St}.) There is no trouble 
``at infinity'', so that the fibration can be extended over
all of $S^3$.  Note that the fibre surface for $\bhat$ 
is the union of $n$ disks with a surface
that is the cover of an annulus branched at $e(\b)$ points, so it 
has Euler characteristic $n-e(\b)$ and hence 
(being connected) genus $g=1-\frac12(n-e(\b)+c)$ if 
$\bhat$ has $c$ components.  This is the same genus 
formula as before when the link of a singularity is considered.

Besides exponent sum there are other representations of $B_n$ 
with applications here.  First recall the \emph{permutation 
representation} $\pi:B_n\to \Ss_n$, which takes 
$\s_j$ to $\transp{j}{j+1}$, $j=1,\dots, n-1$.  
The kernel $\ker\pi$ is the group of \emph{pure braids}; it 
is the fundamental group of the space of \emph{ordered} $n$-tuples 
of distinct complex numbers.  Let $S_n$ be the free abelian 
group of rank $\frac12 n(n-1)$ 
consisting of symmetric $n$-by-$n$ integer matrices 
with $0$ diagonal.  Now, in general, a cycle in $\pi(\b)$ 
corresponds to a component of $\bhat$; and in particular 
the closure of a pure braid consists of $n$ (unknotted) 
components which are naturally ordered $1,\dots, n$.  
Define $\lambda :\ker\pi\to S_n$ by setting 
$\lambda(\b)_{i,j}$ equal to twice the linking number of the 
$i$-th and $j$-th components of $\bhat$, for $\b$ pure.  These 
representations are combined in $\omega:B_n\to S_n\ltimes\Ss_n$,
where in the semidirect product $\Ss_n$ acts on $S_n$ by 
conjugation with the standard permutation matrices, and
$$
\omega(\s_j)=\left([\d_{i,i+1}+\d_{i+1,i}],\transp{i}{i+1}\right),
i=1,\dots,n-1.
$$

Let $\Ss_n$ act diagonally on $\{1,\dots, n\}^2$,  and let 
$|x|\cdot\transp i j$ denote the orbit of (the cyclic subgroup 
generated  by) $x\in\Ss_n$ on $(i,j)$.  Then for
$i\ne j$, $\b\in B_n$, $\omega(\b)=([a_{pq}], x)$,
the sum 
$$
\sum\limits_{(p,q)\in |x|\cdot\transp i j} a_{pq}
$$
is an integer invariant of $\b$, and appropriate sums of such
invariants are conjugacy class invariants.  In particular, when 
$\pi(\b)$ is an $n$-cycle (so that $\b$ is a knot), such a conjugacy 
class invariant arises by summing over pairs $(i,j)$ with a fixed 
constant difference modulo $n$: and this may be seen to be 
precisely twice one of the \emph{self-windings} $\sw_i$ introduced by 
Laufer \cite{Lau}.  Laufer showed that the $\sw_i$ ($i=1,\dots,n$) 
suffice to distinguish the knot types of links of uinibranch 
singularities; in fact, he showed that the Puiseux pairs of a 
branch could be reconstructed from the self-windings.  Simple 
examples show that $\sw=e$ and the $\sw_i$ (and even their slight 
generalizations just given) can't tell apart all quasipositive, 
or even all positive, closed braids.  It is interesting to 
speculate that there might be reasonable representations 
$\lambda_1$ of $\ker\lambda$, $\lambda_2$ of $\ker\lambda_1$, 
\dots, which could somehow be combined into a (faithful?) 
representation of $B_n$ in which quasipositivity might show 
up more clearly than it does in $B_n$ itself.  (Is there any 
relation to Laufer's other numerical link 
invariants \cite{Lau2}? Perhaps $\lambda_1$ can be 
constructed out of linking numbers in branched covers of $S^3$, 
branched over one of the---unknotted!---components of a pure 
braid in which every linking number is $0$, and so on.)

As a final topic, we return to ``knot groups'' of plane curves and 
related matters, from a braid-theoretical point of view.

As before, let $R$ be the compact region of $\C$ bounded by a 
simple closed curve $\g$.  Let $S$ be a compact oriented 
surface-with-boundary.  Then a map $f:S \to R\times\C$, 
or its image $f(S)$, is a \emph{braided surface} of degree 
$n\ge 1$ provided that $\pr_1\circ f:S \to R$ is a 
branched covering of degree $n$: $f$ is a \emph{smooth}, 
\emph{analytic}, or \emph{algebraic} 
braided surface if $f(S)$ is smooth, complex analytic, or 
(complex) algebraic.  Let $V'_f\sub S$ and $V_f\sub R$ denote the 
branch sets of $\pr_1\circ f$, finite sets avoiding $\Bd S$ and 
$\g$; and let $W_f$, $V_f\sub W_f\sub R$ be the set 
$\{z\in R:(\{z\}\times\C)\cap f(S)%
\text{ contains fewer than $n$ points}\}$.  
One can interpret $f ^{-1}$ as a map, as
smooth as $f$, from $R$ into $E_n$.  As remarked earlier, when 
$f^{-1}$ is transverse to $\D$,
then $W_f=V_f$ and $f$ is a smooth braided surface; 
but $f$ can be smooth without
$f$ being transverse to $\D$. 
(Consider non-generic ``vertical'' tangencies.) Nor
need $W_f$ be finite, \emph{but we will always assume that it is,} 
even when $f^{-1}$ is not transverse to $\D$.  
With this proviso, \emph{every braided surface} $f$ is 
a topological (even p.l.) immersion, though not 
necessarily locally flat.  To see this, define the 
\emph{local braid of $f$ at $z\in R$}, denoted $\b_{f, z}\in B_n$, 
to be the homotopy class of the loop 
$\th \mapsto f^{-1}(z +\e e^{i\th})$, 
$0\le 0\le 2\pi$, for any sufficiently small $\e>0$.  
(Since the basepoints of the various copies of $B_n$ 
vary with $z$, $\b_{f,z}$ is really only defined up to
conjugacy.)  This is well-defined when $W_f$ is finite 
(or even as long as $z$ is not an accumulation point 
of $W_f$); of course $\b_{f,z}=1$ if and 
only if $z\in R-W_f$.  For $z\in W_f$, $\bhat_{f ,z}$ 
has strictly fewer than $n$ components, which will be 
grouped into possibly yet fewer unsplittable links.  
Then $f(S)$, above $z$, is embedded in $R\times\C$
like disjoint cones (with distinct vertices) on the 
unsplittable sublinks of $\bhat_{f,z}$.  For
example, if $z\in V_f$ lies under only a simple 
vertical tangent, then $\b_{f,z}$ is a band
(positive or negative), which might as well be 
taken to be $\s_1^{\pm 1}\in B_n$, and $\bhat_{f,z}$ 
is a split link of $n-1$ unknotted components.

Recall (cf. \cite{Bi}) that $B_n$ acts (faithfully) as a group of 
automorphisms of the free group $F_n$ of rank $n$. Explicitly, if 
$F_n = \pi_1(\C\setminus \{ w_1, ..., w_n\}; w_0)$, the acting $B_n$ is 
realized as $\pi_1(E_n-\D;\{w_1,\dots,w_n\})$; on standard free 
generators $x_1,\dots, x_n$ of $F_n$ (positively oriented 
meridians around $w_1,\dots, w_n$), the action~is
$$
x_i\s_i = x_i x_{i+1} x_i^{-1},\thickspace 
x_{i+1}\s_i = x_i,\thickspace
x_j\s_i = x_j \text{ for } j\ne i, i+1.
$$
Pick a basepoint $z_0 \in R\setminus W_f$, and paths from 
$z_0$ to the points $z_1,\dots, z_k$ of $W_f$.   
By these paths, all the local braids can be taken to 
lie in one and the same braid group, namely, 
$\pi_1(E_n\setminus \D; \pr_2((\{z_p\}\times \C)\cap f(S)))$---denote 
by $\bhat'_{f,z}$ these braids.  (Simple vertical tangents, 
for instance, will now give braids $\b'_{f,z}$ which are bands 
that \emph{cannot} all at once be taken to be $\s_1^{\pm1}$.)
It may now be seen that
$$
(x_1,\dots,x_n : x_i \b'_{f,z} = x_i, i = 1,\dots, n, z\in W_f)
$$
is a presentation of the ``knot group'' 
$\pi_1((R \times \C)\setminus f(S);*)$. 
When $f$ is algebraic and $\g$ is a very large circle this is 
really van Kampen's presentation (except for the relation ``at 
infinity'' to which we will return shortly).

A finite presentation of a group, in which each relation sets one 
generator equal to some conjugate of another generator, is a 
\emph{Wirtinger presentation}; a group with a Wirtinger presentation is 
a \emph{Wirtinger group}. Any Wirtinger group has a \emph{simple Wirtinger 
presentation}, in which each relation is of the form 
$x_i x_j x_i^{-1} = x_k$, for not necessarily distinct generators 
$x_i, x_j, x_k$. After possibly adding more generators, and renumbering 
them, one can assume that each relation is of one of the two 
forms $x_i = x_{j+1}$ or $x_i = x_j x_{j+1} x_j^{-1}$, $i <j$.  
These two relations are contributed, respectively, by the action 
on $F_n$ of
$$
(\s_i \s_{i+1}\dotsm\s_{j-1}^\e) \s_j (\s_i \s_{i+1}\dotsm\s_{j-1}^\e)^{-1}, 
\quad\e = {+}1 \text{ or } {-1}.
$$
So every Wirtinger group has a simple Wirtinger presentation 
which is the van Kampen presentation of the fundamental group 
$\pi_1(\{z, w) \in \C^2 : |z|\le 1\} \setminus f(S); *)$ 
for some smooth braided surface $f(S)$ with boundary the 
closure of a quasipositive braid (the product of the bands 
used to achieve the desired relations); and actually $f(S)$ can 
be taken to be non-singular complex analytic.  So we see that \emph{the 
class of knot groups of complex analytic curves in a bidisk is 
exactly the class of Wirtinger groups}, a refinement \cite{Ru2}
of results of Yajima \cite{Ya} and Johnson \cite{Jo} (who weren't 
concerned with complex analytic structures).

If one wishes to investigate knot groups for smooth braided 
surfaces of fixed topological type, one still loses nothing by 
demanding that the surfaces be complex curves: if $f(S)$ is 
smooth, by slight jiggling $f^{-1}$ becomes transverse to $\D$ 
while $f(S)$ moves by an isotopy; then the braids $\b'_{f,z}$  
are all bands, positive or negative; changing all the signs to 
positive reimbeds $S$ as a quasipositive braided surface, and 
therefore, up to isotopy, a complex analytic curve, but it does 
not change the knot group at all, since $x\b^{-1} = x$ is the 
same relation as $x = x\b$.

So far everything has been phrased for braided surfaces over a 
compact (simply-connected) region $R$. If we replace $R$ by all of 
$\C$, much stays the same; it is now appropriate to let $W_f$ be 
infinite, but discrete. It ceases to be clear, however, (at least 
to this author at the present time) that a quasipositive 
``infinite band representation'' can always be realized by an 
entire $n$-valued analytic function. Also, as observed in \cite{Ru1}, 
for compact $R$, at least as far as the boundary closed braid is 
concerned, every $n$-valued analytic function can be assumed to be 
the restriction of an entire $n$-valued \emph{algebraic} function; 
this is certainly not true for $R = \C$, because the ``local braid at 
infinity'' $\b_{f,\infty}$ of an algebraic braided surface over 
$\C$---i.e., the braid over a simple closed curve large enough to 
enclose $V_f$ entirely---is severely restricted.  Its closure, for 
instance, is an iterated torus link (as we saw in the proof of 
the theorem of Abhyankar and Moh, \S 6). And if the projective 
completion of the algebraic braided surface (algebraic curve), in 
$\C\P^2$, meets the line at infinity transversely, one actually 
has $\bhat_{f,\infty}$ the union of $n$ circles of the Hopf 
fibration $S^3 \to \C\P^1$---the braid $\b_{f,\infty}$ is the 
generator of the (infinite cyclic) center of $B_n$ ($n \ge 3$), 
which bears the name $\D^2$ (unfortunately, in this context), 
cf.\ \cite{Bi}. 
Any knot group of a projective plane curve, then, can be 
presented by starting with an expression of $\D^2$ as a product 
$\b(1)\dotsm \b(k)$ in $B_n$, where each $\b(i)$ is conjugate in 
$B_n$ to some local braid associated to the link of a singularity 
(including non-trivial local braids which are associated to the 
unknotted link of a regular point!), then forming the 
presentation
$$
(x_1, \dots, x_n: x_i x_2 \dotsm x_n = 1, x_i\b(j) = x_i, 
i = 1,\dots, n, j =~ 1,\dots, k).
$$
For instance, a quasipositive band representation of $\D^2$ (each 
$\b(i)$ a positive band, that is, conjugate to the nontrivial 
local braid $\s_1$ associated to a simple vertical tangent) 
corresponds to a non-singular curve of degree $n$, and presents 
$\Z/n\Z$. A quasipositive \emph{nodal} band representation, where each 
$\b(i)$ is either a positive band or the square of a positive 
band, corresponds to a node curve; if some $\b(i)$ are cubes of 
positive bands, others squares or first powers, we have a 
\emph{cuspidal} band representation; and so on. There is a mapping from 
the set of strata of $Q_n$ (\S5) into a hierarchy of ``types of 
expressions'' of $\D^2\in B_n$ as products $\b(1)\dotsm\b(k)$; 
Moishezon's problem of normal forms is a first step in the study 
of this mapping, about which little seems to be known. Is it 
onto? An affirmative answer would be a strong generalization of 
Riemann's Existence Theorem. (Again, cf. \cite{Mo}.)

We conclude with three examples. First recall some formulas for 
$\D^2$ in $B_n$ (cf.\ \cite{Bi} or \cite{Mo}): 
$\D^2 = (\s_1\s^2\dotsm\s_{n-1})^n$; also, $\D^2$ is pure, 
and in terms of the standard generators
$$
{A_{ij} = (\s_i\dotsm\s_{j-1})\s_j^2(\s_i\dotsm\s_{j-1})^{-1}, \quad%
1\le i \le j \le n-1,}
$$
of the pure braid group,
$$
\D^2 = A_{1,n-1}A_{1,n-2}\dotsm A_{1,1}A_{2,n-1}%
\dotsm A_{2,2}\dotsm A_{n-1,n-1}.
$$

\emph{Example 1.} Write $\D^2 = \b(1)\dotsm\b(n^2-n)$, 
$\b(i) = \s_{i \mod {n-1}}$, as just given. 
It is easy to see that this expression for $\D^2$ 
does in fact corresponnd to a non-singular curve of 
degree $n$. The corresponding presentation of 
the knot group of the curve includes among its relations
$x_1 x_2\dotsm x_n=1$ and each equality $x_i=x_{i+1}$, 
$i = 1,\dots, n-1$. So the knot group is a quotient of
$\Z/n\Z$; but a simple homological argument
shows that $\Z/n\Z$ is the abelianization of the knot 
group, so the two groups are equal.

\emph{Example 2.} Write $\D^2 = \b(1)\dotsm\b((n^2-n)/2)$, where 
$\b(i)=A_{p,q}$ as above. Each pair $(p, q)$ arises.
The relations in the corresponding presentation say 
that for each pair $p, q$ the generators $x_p, x_{q+1}$ 
commute. (For instance, the action of $A_{1,1} = \s_1^2$ 
on $F_n$ is 
\begin{gather*}
x_1 \s_1^2 = (x_1 x_2 x_1^{-1})\s_1 = x_1 x_2 x_1 x_2^{-1} x_1^{-1}, \\
x_2 \s_1^2 = x_1 \s_1^2 = x_1 x^2 x_1^{-1},\thickspace%
 x_k\s_1^2 = x_k,\thickspace
 k\ne 1,2;
\end{gather*}
and the relations $x_1 = x_1 x_2 x_1 x_2^{-1} x_1^{-1}$ and
$x_2 = x_1 x_2 x_1^{-1}$ both say $x_1$ commutes with $x_2$.) 
The group is free abelian of rank $n-1$. Moishezon sketches a 
proof that this presentation does arise geometrically; another 
proof could be given by the method of \cite{Ru1}. 

\emph{Example 3.} For $n=4$, 
$\D^2 = \s_1 \s_2 \s_3 \s_1 \s_2 \s_3 \s_1 %
\s_2 \s_3 \s_1 \s_2 \s_3$. 
Let us suppress the symbol $\s$, raise subscripts (so 
$k$ denotes $\s_k$), and write, for instance, 
$\SideSet{2 \overline 3}4$ to mean 
$\s_2 \s_3^{-1} \s_4 \s_3 \s_2^{-1}$.
Then, by dogged manipulation, $\D^2 \in B_4$ can be
worked into the form 
$(3\cdot 3\cdot 3) (\SideSet{3 \overline 2}1) 
(1\cdot 1 \cdot 1) (2) (1 \cdot 1 \cdot 1) (\SideSet{32}1)$, 
the product of three positive bands and three ``cusps'' (cubes 
of positive bands). The corresponding presentation, before adjoining 
the relation at infinity, presents the group of the $5$-twist spun 
trefoil (as has been remarked by Dewitt Sumners); with 
that relation, $x_1 x_2 x_3 x_4 = 1$, the group becomes the 
non-abelian group of order 12, $(a, b: aba = bab, a^4 = 1, a^2=b^2)$. 
This is the correct group \cite{Z} for a tricuspidal cubic curve,
and presumably the given ``quasipositive cuspidal band 
representation'' really arises geometrically, but I have 
not had the courage to check this.---Similarly, for
$n = 6$, $\D = 123451234123121$, which can be written as
$(1\cdot 1 \cdot 1) 
(\SideSet{{\overline 1}2}1)  
(3\cdot 3 \cdot 3)
(\SideSet{{\overline 3}4}3)
(5\cdot 5\cdot 5) 
(\SideSet{\overline{133}}2)
(\SideSet{\overline{355}}4)
(\SideSet 2 3)
(\SideSet 4 5)$;
the presentation for 
the square of this, with the relation at infinity, is at an 
intermediate stage 
\begin{multline*}
(x_1, x_2, x_3, x_4, x_5, x_6 : %
x_1 = x_3 = x_5, 
x_2 = x_4 = x_6,\\
x_1 x_2 x_1 = x_2 x_1 x_2, 
x_1 x_2 x_3 x_4 x_5 x_6 = 1)
\end{multline*}
which becomes 
$(a,b: a^2=b^3=1)$, the group given in \cite{Z} 
for a sextic with six cusps on a conic. 
On the other hand, a less symmetrical way to write 
$\D^2 \in B_6$ is as 
\begin{multline*}
(\SideSet{2{\overline 1}2}3)
(4)(5)   
(2\cdot 2 \cdot 2) 
(\SideSet 1{2}\cdot \SideSet 1{2} \cdot \SideSet 1{2}) 
(\SideSet 3 2)
(\SideSet{431}2)
(1\cdot 1 \cdot 1)
(\SideSet{4{\overline 3}2}1)
(\SideSet{44}5)\\
(\SideSet{442}3)
(4\cdot 4 \cdot 4)
(1\cdot 1 \cdot 1) 
(\SideSet{2{\overline 1}2}3)
(\SideSet{2{\overline 1}2}3)
(1\cdot 1 \cdot 1) 
(\SideSet{22}1)
(2),
\end{multline*}
which presents the abelian group $\Z/6\Z$ which 
\cite{Z} gives for a sextic with six cusps \emph{not} all on the same conic.

\renewcommand{\bysame}%
{\leavevmode\raisebox{1.5pt}{\hbox to1.8em{\hrulefill}}\medspace\medspace}

\medskip
\obeylines
{\small\noindent Lee Rudolph
\smallskip
\noindent\phantom{Lee }POB 251
\noindent\phantom{Lee }Adamsville, R.I.\ 02801 (USA)}


\begin{thebibliography}{Lau 2x}

\bibitem[Abe]{Abe}
{\sc Abelson}, Harold. {Fundamental groups of plane curves and their duals}.
  \emph{Indiana Univ. Math. J. 25} (1976), No.~1, 65--67.

\bibitem[Ab]{Ab}
{\sc Abhyankar}, Shreeram.  {Tame coverings and fundamental groups of algebraic
  varieties. {V}. {T}hree cuspidal plane quartics}. 
  \emph{Amer. J. Math. 82} (1960), 365--373.

\bibitem[A-M]{A-M}
{\sc Abhyankar}, S. and Tzuong~Tsieng {\sc Moh}. 
  {Embeddings of the line in the plane}. 
  \emph{J. reine ange. Math. 276} (1975), 148--166.

\bibitem[A-S]{A-S}
{\sc Abhyankar}, S. and Balwant {\sc Singh}. 
  {Embeddings of certain curves in the affine plane}. 
  \emph{Amer. J. Math. 100} (1978), 99--175.

\bibitem[Bi]{Bi}
{\sc Birman}, Joan. \emph{Braids, Links, and Mapping Class Groups}. 
  Ann. Math. Studies 82 (1975), Princeton University Press.

\bibitem[Bi-W]{Bi-W}
{\sc Birman}, J. and R.~F. {\sc Williams}. 
  {Knotted periodic orbits in dynamical systems. I: Lorenz's 
  equations}. \emph{Topology 22} (1983), 47--82.

\bibitem[Bl]{Bl}
{\sc Bliss}, Gilbert~Ames. \emph{Algebraic Functions}.  Amer. Math. Soc.
  Colloquium Publications, vol. XVI (1933).

\bibitem[Cha]{Cha}
{\sc Chang}, Haichau. {On two classical facts by Zariski and 
  van Kampen}. \emph{Chinese J. Math. 7} (1979), 153--161.

\bibitem[Che]{Che}
{\sc Cheniot}, D. {Une d\'emonstration du th\'eor\`eme de {Z}ariski sur les
  sections hyperplanes d'une hypersurface projective et du th\'eor\`eme de
  {V}an {K}ampen sur le groupe fondamental du compl\'ementaire d'une courbe
  projective plane}. \emph{Compositio Math. 27} (1973), 141--158.

\bibitem[De]{De}
{\sc Deligne}, P. {Le groupe fondamental du compl\'ement d'une courbe plane
  n'ayant que des points doubles ordinaires est ab\'elien.}
  \emph{S{\'e}m. Bourbaki}, Nov. 1979.

\bibitem[E-N]{E-N-preprint}
{\sc Eisenbud}, D.\ and Walter {\sc Neumann}. 
  {Fibering {I}terated {T}orus {L}inks},
  \emph{To appear.}

\bibitem[F-H]{F-H}
{\sc Fulton}, W.\ and J.\ {\sc Hansen}. 
  {A connectedness theorem for projective
  varieties, with applications to intersections 
  and singularities of mappings}.
  \emph{Ann. Math. 110} (1979), 159--166.

\bibitem[G-R]{G-R}
{\sc Gunning}, R.\ C.\ and Hugo {\sc Rossi}. 
  \emph{Analytic Functions of Several Complex
  Variables}. Prentice-Hall, 1965.

\bibitem[H]{H}
  {\sc Hitotumatu}, Sin. Some recent topics in 
  several complex variables by the Japanese school. 
  (Report of work by T.\ Nishino), in 
  \emph{Proceedings of the
  Romanian-Finnish Seminar on Teichm\"uller Spaces and 
  Quasi-Conformal Mappings}, Brasov, Romania, 1969, 
  Acad. Soc. Rep. Romania (1971), 187--191.

\bibitem[Jo]{Jo}
{\sc Johnson}, D. \emph{Homomorphs of knot groups}.  Preprint,
  Jet Propulsion Laboratory, Pasadena, California (1978);
  \emph{Proc. Amer. Math. Soc. 78} (1980), 135--138.

\bibitem[Lau]{Lau}
{\sc Laufer}, Henry B. {On the number of singularities of an analytic curve}.
  \emph{Trans. Amer. Math. Soc. 136} (1969), 527--535.

\bibitem[Lau 2]{Lau2}
\bysame {Some numerical link invariants}. \emph{Topology 10}
  (1971), 119--131.

\bibitem[L\^e]{Le}
{\sc L\^e~Dung Trang}, {Sur les n{\oe}uds alg\'ebriques}.
  \emph{Compositio Math.\ 25} (1972), 281--321.

\bibitem[Li]{Li}
{\sc Litherland}, R.~A. {Signatures of iterated torus knots}.  In
  \emph{Topology of low-dimensional manifolds},
  Lecture Notes in Mathematics 772 (1979), 
  Springer, Berlin, 71--84.

\bibitem[Mi 1]{Mi1}
{\sc Milnor}, J. W. \emph{Morse Theory}.  Ann. Math. Studies 51 (1965), 
  Princeton University Press.

\bibitem[Mi 2]{Mi2}
\bysame \emph{Singular Points of Complex Hypersurfaces}. 
  Ann. Math. Studies 61 (1969), Princeton University Press.

\bibitem[Mo]{Mo}
{\sc Moishezon}, B.~G. {Stable branch curves and braid monodromies}, 
  In \emph{Algebraic Geometry (Chicago, Illinois, 1980)}, 
  Lecture Notes in Mathematics 862 (1981), 107--192, Springer, Berlin.

\bibitem[Ru 1]{Ru1}
{\sc Rudolph}, Lee. {Algebraic functions and closed braids}. 
  \emph{Topology 22} (1983), 191--202.

\bibitem[Ru 2]{Ru2}
\bysame {Braided surfaces and {S}eifert ribbons for closed braids}.
  \emph{Comment. Math. Helv. 58} (1983), no.~1, 1--37.

\bibitem[Ru 3]{Ru3}
\bysame {Constructions of quasipositive knots and links, {I}}. 
  In \emph{Knots, Braids and Singularities (Plans-sur-Bex, 1982)}, 
  Univ. Gen\`eve, Geneva, 1983, 233--245.

\bibitem[Ru 4]{Ru4}
\bysame {Embeddings of the line in the plane}. 
  \emph{J. reine ange. Math. 337} (1982), 113--118.

\bibitem[Ru 5]{Ru5}
\bysame {Non-trivial positive braids have positive signature}.
  \emph{Topology 21} (1982), no.~3, 325--327.

\bibitem[Ru 6]{Ru6}
\bysame {Question printed in ``Queries'' column}. 
   \emph{Notices Amer. Math.  Soc. 23} (1976), p.\ 410 
   (from a problem list compiled at the Special Session 
   on Knot Theory, 1976 Summer Meeting of A.M.S., Toronto).

\bibitem[St]{St}
{\sc Stallings}, J. {Constructions of fibred knots and links}.  In
  \emph{Algebraic and Geometric Topology} (Proc. Sympos. 
  Pure Math.\ XXXII) part 2 (1978), 55--60, Amer. Math. Soc., 
  Providence, R.\ I.

\bibitem[vK]{vK}
{\sc van Kampen}, E.~R. {On the fundamental group of an 
  algebraic curve}. \emph{Amer. J. Math 55} (1933), 255--260.

\bibitem[Ya]{Ya}
{\sc Yajima}, T. {Wirtinger presentations of knot groups}.
  \emph{Proc.\ Japan Acad.\ Sci.\ 46} (1970), 997--1000.

\bibitem[Z]{Z}
{\sc Zariski}, O. {On the problem of existence of algebraic functions of two
  variables possessing a given branch curve}. 
  \emph{Amer. J. Math 51} (1929), 305--328.
\smallskip

{\hfill \emph{(\,Re{\c c}u le 30 september 1982\,)}}

\bigskip
\centerline{ADDITIONAL REFERENCES}
\smallskip

\bibitem[B-O]{B-O}
{\sc Boileau}, Michel and Stepan Yu. {\sc Orevkov}.
Quasipositivit\'e d'une courbe analytique dans une boule 
pseudo-convexe.  \emph{C. R. Acad. Sci. Paris 332} (2001), 825--830. 

\bibitem[E-N 2]{E-N-published}
{\sc Eisenbud}, D.\ and Walter {\sc Neumann}.
\emph{Three-dimensional Link Theory and Invariants of Plane 
Curve Singularities}. Ann.\ Math.\ Studies 110, 1985

\bibitem[Ha]{Ha}
{\sc Harris}, Joe. On the Severi problem.
\emph{Invent. Math. 84} (1986), 445--461.

\bibitem[K-M]{K-M}
{\sc Kronheimer}, P.B.\ and T.S.\ {\sc Mrowka}.
Gauge theory for embedded surfaces. I.
\emph{Topology 32} (1993), 773--826.

\bibitem[N]{N}
{\sc Neumann}, Walter D.
Complex algebraic plane curves via their links at infinity.
\emph{Invent. Math. 98} (1989), 445--489.

\bibitem[N-R]{N-R}
{\sc Neumann}, Walter and Lee {\sc Rudolph}.
Unfoldings in knot theory.
\emph{Math. Ann. 278} (1987), 409--439

\bibitem[O]{O}
{\sc Orevkov}, Stepan Yu. 
The fundamental group of the complement of a plane algebraic curve.
\emph{Mat. Sbornik 137} (1988), 260--270.

\bibitem[Ru 7]{Ru7}
{\sc Rudolph}, Lee. Quasipositivity as an obstruction to sliceness.
  \emph{Bull. Amer. Math. Soc. 29} (1993), 51--59.

\bibitem[Suz]{Suz}
{\sc Suzuki}, Masakazu.
{Propri\'et\'es topologiques des polyn\^omes de deux variables
complexes, et automorphismes alg\'ebriques de l'espace ${\bf {C}}\sp{2}$}.
\emph{J. Math. Soc. Japan 26} (1974), 241--257.

\bibitem[Z-L]{Z-L}
{\sc Za{\u\i}denberg}, M.\ G.\ and V.\ Ya.\ {\sc Lin}.
An irreducible, simply connected algebraic curve in
${\bf {C}}\sp{2}$\ is equivalent to a quasihomogeneous curve.
\emph{Dokl. Akad. Nauk SSSR 271} (1983), 1048--1052.

\end{thebibliography}
\end{document}